\documentclass[11pt]{amsart}
\usepackage{a4wide}
\usepackage{amsmath, amsfonts, amssymb,graphicx}
\usepackage{mathrsfs,url}

\date{\today}

\author{Manuel Bodirsky}
    \address{Laboratoire d'Informatique  (LIX), CNRS UMR 7161\\
    \'{E}cole Polytechnique \\91128 Palaiseau\\
    France}
    \email{bodirsky@lix.polytechnique.fr}
    \urladdr{http://www.lix.polytechnique.fr/~bodirsky/}
\author{Michael Pinsker}
    \address{\'{E}quipe de Logique Math\'{e}matique\\ Universit\'{e} Diderot - Paris 7\\
	UFR de Math\'{e}matiques\\
	75205 Paris Cedex 13, France}
    \email{marula@gmx.at}
    \urladdr{http://dmg.tuwien.ac.at/pinsker/}
    \thanks{The first author has received funding from the European Research Council under the European Community's Seventh Framework Programme (FP7/2007-2013 Grant Agreement no. 257039). The second author is grateful for support through an APART-fellowship of the Austrian Academy of Sciences.}

\title[Topological Birkhoff]{Topological Birkhoff}


\DeclareMathOperator{\fin}{fin}

\DeclareMathOperator{\pol}{Pol}
\DeclareMathOperator{\Emb}{Emb}
\DeclareMathOperator{\Clo}{Clo}
\DeclareMathOperator{\HSP}{HSP}
\DeclareMathOperator{\HHH}{H}
\DeclareMathOperator{\PPP}{P}
\DeclareMathOperator{\PPPfin}{P^{\fin}}
\DeclareMathOperator{\SSS}{S}

\newcommand{\HSPfin}{\HSP^{\fin}}
\newcommand{\rest}{{\upharpoonright}}

\newcommand{\cC}{\ensuremath{\mathcal{C}}}
\newcommand{\cV}{\ensuremath{\mathcal{V}}}

\DeclareMathOperator{\Betw}{Betw}

\DeclareMathOperator{\OIT}{1IN3}

\newcommand{\ignore}[1]{}

\newcommand{\To}{\rightarrow}
\newcommand{\nin}{\notin}
\newcommand{\inv}{^{-1}}
\newcommand{\mult}{\times}
\newcommand{\uk}{^{(k)}}
\newcommand{\cloa}{\Clo(\A)}
\newcommand{\clob}{\Clo(\B)}
\newcommand{\ccloa}{\overline{\Clo(\A)}}
\newcommand{\cclob}{\overline{\Clo(\B)}}

\DeclareMathOperator{\Csp}{CSP}

\DeclareMathOperator{\id}{id}

\DeclareMathOperator{\Aut}{Aut}
\DeclareMathOperator{\End}{End}
\DeclareMathOperator{\Pol}{Pol}

\newcommand{\A}{{\bf A}}
\newcommand{\B}{{\bf B}}
\newcommand{\C}{{\bf C}}
\newcommand{\SSSS}{{\bf S}}
\newcommand{\CC}{{\mathscr C}}
\newcommand{\D}{\mathscr D}

\newcommand{\F}{\mathscr F}
\newcommand{\G}{\mathscr G}
\renewcommand{\H}{\mathscr H}

\renewcommand{\O}{\mathscr O}

\theoremstyle{plain}

    \newtheorem{thm}{Theorem}
    \newtheorem{theorem}[thm]{Theorem}

    \newtheorem{lem}[thm]{Lemma}

    \newtheorem{prop}[thm]{Proposition}
    \newtheorem{proposition}[thm]{Proposition}

    \newtheorem{cor}[thm]{Corollary}
    \newtheorem{corollary}[thm]{Corollary}

\theoremstyle{definition}

    \newtheorem{definition}[thm]{Definition}

    \newtheorem{nota}[thm]{Notation}

\begin{document}
\begin{abstract}
	One of the most fundamental mathematical contributions of Garrett Birkhoff is the HSP theorem, which implies that 
	a finite algebra $\B$ satisfies all equations that hold in a finite algebra $\A$
	of the same signature 
	if and only if $\B$ is a homomorphic image of a subalgebra of a finite power of $\A$. On the other hand, if $\A$ is infinite, then in general one needs to take an \emph{infinite} power in order to obtain a representation of $\B$ in terms of $\A$, even if $\B$ is finite.
	
	We show that by considering the natural topology on the functions of $\A$ and $\B$ in addition to the equations that hold between them, one can do with finite powers even for many interesting infinite algebras $\A$. More precisely, we prove that if $\A$ and $\B$ are at most countable algebras which are oligomorphic, then the mapping which sends each function from $\A$ to the corresponding function in $\B$ preserves equations and is \emph{continuous} if and only if $\B$ is a homomorphic image of a subalgebra of a \emph{finite} power of $\A$.
	
	Our result has the following consequences in model theory and in theoretical computer science: two $\omega$-categorical structures are primitive positive bi-interpretable if and only if their topological polymorphism clones are isomorphic. In particular, the complexity of the constraint satisfaction problem of an $\omega$-categorical structure only depends on its topological polymorphism clone.
	\end{abstract}
\maketitle
\section{Introduction}
The algebraic result we present has a motivating application 
in model theory, which in turn has implications
for the study of the computational complexity of constraint satisfaction problems in theoretical computer science. We start our introduction with
this model-theoretic perspective on our result, and describe the central algebraic theorem of this article later in the introduction, in Section~\ref{sect:ua-intro}.

\subsection{The model-theoretic perspective}
A countable structure $\Gamma$ is called 
\emph{$\omega$-categorical} iff all countable models of the first-order theory of $\Gamma$ are isomorphic to $\Gamma$. A substantial amount of information about an $\omega$-categorical structure $\Gamma$ is already coded into the automorphism group $\Aut(\Gamma)$ of $\Gamma$, 
viewed abstractly as a topological group whose topology is the topology of pointwise convergence. In particular, 
Ahlbrandt and Ziegler~\cite{AhlbrandtZiegler} proved that two countable $\omega$-categorical structures are \emph{first-order bi-interpretable} if and only if their automorphism groups are isomorphic as topological groups. The concept of interpretation we use here is standard~\cite{HodgesLong}, and will be recalled in Section~\ref{sect:interpretations}. 

Recently, the following variant of the theorem of Ahlbrandt and Ziegler has been shown, replacing the automorphism group by the endomorphism monoid (which, of course, contains more information about the original structure than the automorphism group)~\cite{BodJunker}:
 two $\omega$-categorical structures $\Gamma$ and $\Delta$ without constant endomorphisms are \emph{existential positive bi-interpretable} (i.e., bi-interpretable by means of existential positive first-order formulas) if and only if their endomorphism monoids $\End(\Gamma)$ and $\End(\Delta)$ are isomorphic as abstract topological monoids, i.e., iff there exists a bijective function $\xi\colon \End(\Gamma)\To \End(\Delta)$ which sends the identity function on $\Gamma$ to the identity on $\Delta$, which satisfies $\xi(f\circ g)=\xi(f)\circ\xi(g)$ for all $f,g\in\End(\Gamma)$, and such that both $\xi$ and its inverse are continuous.
 
In the same paper, it is stated as an open problem whether this statement can be modified further to characterize \emph{primitive positive bi-interpretability}, if one replaces the endomorphism monoid by the \emph{polymorphism clone}. 
A primitive positive interpretation is a first-order interpretation where all the involved formulas are \emph{primitive positive}, i.e., of the form $\exists x_1,\dots,x_n (\phi_1 \wedge \dots \wedge \phi_m)$ where $\phi_1,\dots,\phi_m$ are atomic formulas. A \emph{polymorphism} of a structure $\Gamma$ is a homomorphism from $\Gamma^k$ to $\Gamma$ for some finite $k\geq 1$; the polymorphism clone $\Pol(\Gamma)$ of $\Gamma$ is the set of all polymorphisms of $\Gamma$  and contains, in particular, at least the information of $\End(\Gamma)$, which is the unary fragment of $\Pol(\Gamma)$. In general, a \emph{(concrete) clone} is a set of finitary functions on a fixed set which contains all projections and which is closed under composition; it is not hard to see that $\Pol(\Gamma)$ is a clone in this sense. Moreover, $\Pol(\Gamma)$ is a closed subset of the topological space $\O_\Gamma=\bigcup_{k\geq 1} \Gamma^{\Gamma^k}$ of all finitary functions on $\Gamma$, just like $\Aut(\Gamma)$ is a closed subset of the space of all permutations on $\Gamma$ and $\End(\Gamma)$ is a closed subset of the space of unary functions on $\Gamma$. The topology of $\O_\Gamma$ is obtained by viewing this space as the sum space of the spaces $\Gamma^{\Gamma^k}$, and each $\Gamma^{\Gamma^k}$ as a power of $\Gamma$, which itself is taken to be discrete. Similarly to automorphism groups and endomorphism monoids, where we distinguish between the concrete permutation groups and transformation monoids on the one hand and abstract topological groups and topological monoids with their laws of composition and topology on the other hand,
 polymorphism clones can be viewed abstractly as \emph{topological clones} carrying an algebraic and a topological structure. The topology on $\Pol(\Gamma)$ is inherited from the space $\O_\Gamma$; note that each $\Gamma^{\Gamma^k}$, and in fact also $\O_\Gamma$, is homeomorphic to the Baire space, and that therefore the space
 $\Pol(\Gamma)$ is a closed subspace of the Baire space. The algebraic structure of $\Pol(\Gamma)$ is that of a multi-sorted algebra with operations that correspond to the composition of the elements of $\Pol(\Gamma)$ and constant symbols corresponding to the projections. We can avoid a formal description of this ghastly structure here (and refer the interested reader, for example, to~\cite{AbstractClones} or the survey paper \cite{GoldsternPinsker}) since we only need the very natural notion of a \emph{homomorphism between clones $\CC,\D$}: these are  functions $\xi\colon\CC\To\D$ which send every  projection in $\CC$ to the corresponding projection in $\D$, and such that $\xi(f(g_1,\ldots,g_n))=\xi(f)(\xi(g_1),\ldots,\xi(g_n))$ for all $n$-ary $f\in \CC$ and all $m$-ary $g_1,\ldots g_n\in\CC$. In particular, two polymorphism clones $\Pol(\Gamma),\Pol(\Delta)$ are isomorphic as topological clones iff there exists a bijection $\xi$ from $\Pol(\Gamma)$ onto $\Pol(\Delta)$ such that both $\xi$ and its inverse are continuous clone homomorphisms. The above-mentioned problem in~\cite{BodJunker} asked whether for two $\omega$-categorical structures $\Gamma,\Delta$ having isomorphic polymorphism clones and being primitive positive bi-interpretable is one and the same thing.

Besides the theoretical interest they might have, primitive positive interpretations are additionally motivated
by an application in theoretical computer science: every relational structure $\Gamma$ with a finite signature defines a computational problem, called the \emph{constraint satisfaction problem of $\Gamma$} and denoted by $\Csp(\Gamma)$, and it is known that when a relational structure $\Delta$ has a primitive positive interpretation in a relational structure $\Gamma$, then $\Csp(\Delta)$ has a polynomial-time reduction
to $\Csp(\Gamma)$. Very general and deep complexity classification results rely on this fact -- see for example the collection of survey articles in~\cite{CSPSurveys}; 
more on this application can be found in Section~\ref{sect:csp}.

In this paper, we give an affirmative answer to the question from~\cite{BodJunker} about primitive positive interpretability.  
A \emph{reduct} of a structure $\Delta'$ is a structure on the same domain obtained by forgetting some relations or functions of $\Delta'$. We  prove the following.

\begin{theorem}\label{thm:main}
Let $\Gamma$ and $\Delta$ be finite or countable $\omega$-categorical structures. Then: 
\begin{itemize} 
\item $\Delta$ has a primitive positive interpretation in 
$\Gamma$ if and only if 
$\Delta$ is a reduct of a finite or $\omega$-categorical structure $\Delta'$ 
such that there exists a continuous homomorphism from $\Pol(\Gamma)$ into $\Pol(\Delta')$ whose image is dense in $\Pol(\Delta')$. 
\item $\Gamma$ and $\Delta$ are primitive positive bi-interpretable if and only if their polymorphism clones are isomorphic as topological clones.
\end{itemize} 
\end{theorem}
It follows from this theorem and the remarks above that the computational complexity of the constraint satisfaction problem for a relational structure in a finite language only depends on its
topological polymorphism clone.

\begin{corollary}\label{cor:csp-topo-clone}
	Let $\Gamma$ and $\Delta$ be finite or countable $\omega$-categorical relational structures with finite signatures. If $\Pol(\Gamma)$ and $\Pol(\Delta)$ are isomorphic as topological clones, then $\Csp(\Gamma)$ and $\Csp(\Delta)$ are polynomial-time equivalent.
\end{corollary}

\subsection{Topological Birkhoff}
\label{sect:ua-intro}
To prove Theorem~\ref{thm:main} we show an algebraic result
which is of independent interest 
and which can be seen as a topological version of Birkhoff's HSP theorem.

An \emph{algebra} is a structure with a purely functional signature. The \emph{clone of an algebra $\A$} with signature $\tau$, denoted by $\Clo(\A)$, is the set of all functions with finite arity on the domain $A$ of $\A$ which can be written as $\tau$-terms over $\A$. More precisely, every abstract $\tau$-term $t$ induces a function $t^\A$ on $A$, and $\Clo(\A)$ consists precisely of the functions of this form.

Let $\A$, $\B$ be algebras of the same signature $\tau$. The assignment $\xi$ from $\Clo(\A)$ to $\Clo(\B)$ which sends every element $t^\A$ of $\Clo(\A)$ to $t^\B$ is a well-defined function if and only if for all $\tau$-terms $s,t$ we have that $s^\B=t^\B$ whenever $s^\A=t^\A$. In that case, it is in fact a surjective homomorphism between clones; we then call $\xi$ the \emph{natural homomorphism} from $\Clo(\A)$ onto $\Clo(\B)$.

When $\cC$ is a class of algebras with common signature $\tau$, then $\PPP(\cC)$ denotes the class of all products of algebras from $\cC$, $\PPPfin(\cC)$ denotes the class of all \emph{finite} products of algebras from $\cC$,
$\SSS(\cC)$ denotes the class of all subalgebras of algebras from $\cC$, and $\HHH(\cC)$ denotes the class of all 
homomorphic images of algebras from $\cC$. A \emph{pseudovariety} is a class $\cV$ of algebras
of the same signature 
such that $\cV = \HHH(\cV) = \SSS(\cV) = \PPPfin(\cV)$, i.e., a class closed under homomorphic images, subalgebras, and finite products; the pseudovariety \emph{generated} by a class of algebras $\cC$ (or by a single algebra $\A$) is the smallest pseudovariety that contains $\cC$ (contains $\A$, respectively). 
For \emph{finite} algebras, Birkhoff's HSP theorem takes the following form (see Exercise 11.5 in combination with the proof of Lemma 11.8 in~\cite{BS}).

\begin{theorem}[Birkhoff]\label{thm:birkhoff}
Let $\A,\B$ be finite algebras with the same signature. 
Then the following three statements are equivalent. 
\begin{enumerate}
\item The natural homomorphism from $\Clo(\A)$ onto $\Clo(\B)$ exists. 
\item $\B \in \HSPfin(\A)$.
\item $\B$ is contained in the pseudovariety generated by $\A$.
\end{enumerate}
\end{theorem}
When $\A$ and $\B$ are of arbitrary cardinality, then the equivalence of $(2)$ and $(3)$ still holds; however, if one wants to maintain equivalence with item~(1), then another version of Birkhoff's theorem states that one has to replace
finite powers by arbitrary powers in the second item,
that is, one has to replace $\HSPfin(\A)$
by $\HSP(\A)$; the third item has to be adapted using the notion of a \emph{variety} of algebras, i.e., a class of algebras of common signature 
closed under the operators $\HHH$, $\SSS$ and $\PPP$.

Our topological variant of Birkhoff's theorem shows
that one can keep finite powers for a large class of infinite algebras if one additionally
requires that the natural homomorphism from $\Clo(\A)$ onto $\Clo(\B)$
is continuous when we view $\Clo(\A)$ and $\Clo(\B)$ as 
topological clones as described above.
 
 A permutation group $\G$ on a countable set $A$ is called \emph{oligomorphic} iff for each finite $n\geq 1$, the componentwise action of $\G$ on $A^n$ has finitely many orbits. In our context it is worth noting that the theorem of Ahlbrandt and Ziegler implies that being oligomorphic is a property of the abstract topological group $\G$, i.e., for isomorphic permutation groups $\G$ and $\H$ one is oligomorphic iff the other one is; for a characterization by abstract properties see~\cite{Tsankov}.
An algebra $\A$ is called \emph{oligomorphic} iff 
 the unary invertible operations in $\Clo(\A)$,
 that is, the unary bijective operations whose inverse is also in $\Clo(\A)$, form an oligomorphic permutation group. We call it \emph{locally oligomorphic} iff the topological closure $\overline{\Clo(\A)}$ in the space $\O_\A$ of all finitary functions on the domain of $\A$ is oligomorphic. Clearly, oligomorphic algebras are also locally oligomorphic; the algebra on a countable set $A$ which has all unary operations on $A$ which are not permutations is an example which shows that the two notions are not equivalent.

One of the motivations for oligomorphic groups is the
theorem of Engeler, Svenonius, and Ryll-Nardzewski (see e.g.~the textbook~\cite{HodgesLong}): the automorphism group of a countable structure $\Gamma$ is oligomorphic if and only if $\Gamma$ is $\omega$-categorical. 
This implies that any \emph{polymorphism algebra} of $\Gamma$, i.e., any algebra on the domain on $\Gamma$ whose functions are precisely the elements of $\Pol(\Gamma)$ indexed in some arbitrary way, is oligomorphic if and only if $\Gamma$ is $\omega$-categorical; note that such polymorphism algebras are oligomorphic if and only if they are locally oligomorphic, since their clone $\Pol(\Gamma)$ is always a closed subset of $\O_\Gamma$. It is not hard to see that all algebras in the pseudovariety generated by an oligomorphic (locally oligomorphic) algebra are again oligomorphic (locally oligomorphic). In this paper, we will prove the equivalence of (1) and (2) in the following theorem, which is a topological characterization of pseudovarieties of oligomorphic algebras. As mentioned above, the equivalence of (2) and (3) holds for arbitrary algebras $\A, \B$ and is well-known from Birkhoff's work.

 \begin{theorem}\label{thm:topo-birk}
 Let $\A,\B$ be locally oligomorphic or finite algebras with the same signature. 
Then the following three statements are equivalent. 
\begin{enumerate}
\item The natural homomorphism from $\Clo(\A)$ onto $\Clo(\B)$ exists and is continuous. 
\item $\B \in \HSPfin(\A)$.
\item $\B$ is contained in the pseudovariety generated by $\A$.
\end{enumerate}
\end{theorem}


\begin{figure}[h]
\begin{center}
\includegraphics[scale=.5]{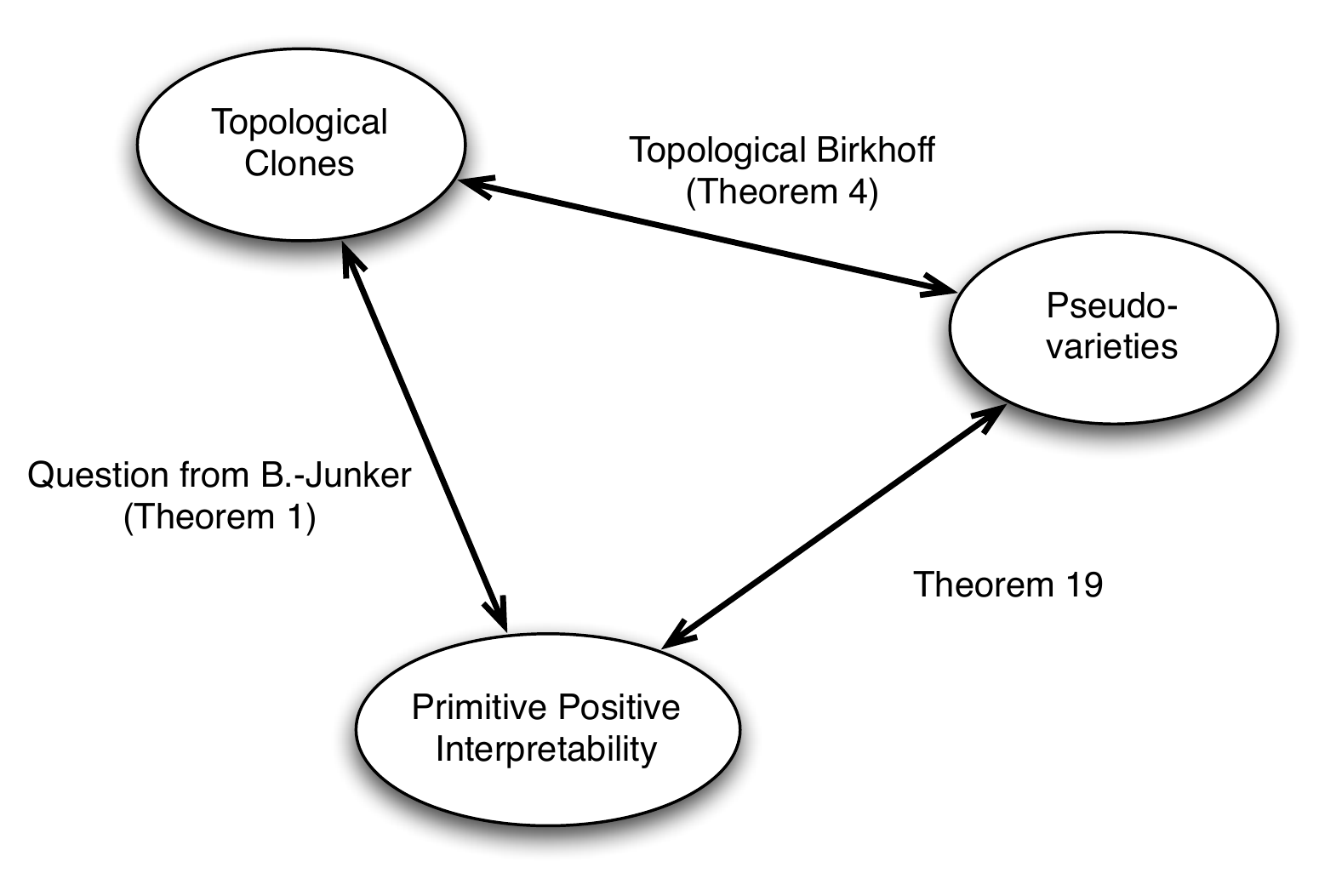} 
\end{center}
\caption{Topological clones, primitive positive interpretations, and pseudovarieties.}
\label{fig:diagram}
\end{figure}



Note that Theorem~\ref{thm:birkhoff} really is a special case of Theorem~\ref{thm:topo-birk}, since the topology of any clone on a finite set is discrete, and hence the natural homomorphism from the clone of a finite algebra to that of another algebra is always continuous.

We will see in Section~\ref{sect:interpretations}
how to derive Theorem~\ref{thm:main} from Theorem~\ref{thm:topo-birk} and a certain correspondence between primitive positive interpretations and pseudovarieties -- confer also Figure~\ref{fig:diagram}. 

\subsection{Related Work}
Pseudovarieties consisting of \emph{finite} algebras have been studied by 
many researchers in many different contexts, and are important in particular in formal language theory.
There is also an equational characterization for pseudovarieties of finite algebras, the Eilenberg-Sch\"utzenberger theorem~\cite{EilenbergSchuetzenberger}.
The topology used in subsequent publications~\cite{Banaschewski,Reitermann} concerning pseudovarieties of finite algebras is different from the topology that we use here; also note that our results are about pseudovarieties 
that also contain infinite algebras.

\cite{MasulovicPech}, in connection with pioneering work on homomorphism-homogeneous structures~\cite{CameronNesetril}, introduce the notion of \emph{weakly oligomorphic} for relational structures via their endomorphism monoid. It could make sense to use their notion for algebras as another weakening of ``oligomorphic'' rather than ``locally oligomorphic''. In fact, every weakly oligomorphic algebra would then be locally oligomorphic, and yet it is quite possible that Theorem~\ref{thm:topo-birk} holds and can be proven our same methods even for the class of weakly oligomorphic algebras. However, we gave preference to ``locally oligomorphic'' with its more group theoretic flavor.

In our proof of Theorem~\ref{thm:topo-birk} we will work with the closure $\overline{\Clo(\A)}$ of $\Clo(\A)$ in $\O_\A$ rather than with $\Clo(\A)$ itself, allowing for a certain compactness argument. Even when the functions of $\A$ are assumed to form, say, a closed set, $\Clo(\A)$ can be topologically complicated: \cite{GPS} give an example of a (topologically) closed algebra $\A$ which has only unary operations and for which $\Clo(\A)$ is not a Borel set.

\subsection{Outline of the paper}
This introduction is followed by Section~\ref{sect:mainproof}, in which we will provide the proof of Theorem~\ref{thm:topo-birk}. We then give some examples in Section~\ref{sect:examples1} which examine the differences between continuous and non-continuous clone homomorphisms in our context. Section~\ref{sect:interpretations} brings us back to the model-theoretic perspective in more detail and links Theorems~\ref{thm:main} and  \ref{thm:topo-birk}. We will provide concrete instances of Theorem~\ref{thm:main} in Section~\ref{sect:examples2}. In Section~\ref{sect:csp} we discuss applications to constraint satisfaction problems; the discussion will be followed by a concrete example in Section~\ref{sect:examples3}. We conclude the paper with an outlook and open problems in Section~\ref{sect:conclusion}.

\subsection{Further conventions}

All $\omega$-categorical structures in this paper are assumed to be countable.

If $\F$ is a set of finitary functions on a set and $k\geq 1$, then we write $\F^{(k)}$ for the $k$-ary functions in $\F$; this applies in particular to $\Pol(\Gamma)$ and $\Clo(\A)$.

For an $n$-tuple $a$ and $1\leq i\leq n$, we write $a_i$ for the $i$-th component of $a$. We do not always distinguish between the domain of a structure and the structure itself, so we write things like ``$a\in\Gamma$'' to refer to an element of $\Gamma$. In the case of algebras, however, we also write $A$ for the domain of $\A$.

When we write $fg$ for a composite of unary functions $f,g$, we mean that $g$ is applied first.


\section{Pseudovarieties and Topological Clones}\label{sect:mainproof}

\subsection{Continuity of the natural homomorphism}

The following lemma shows the easy direction of the equivalence of Theorem~\ref{thm:topo-birk}, namely that (2) implies (1).

 \begin{prop}\label{prop:hin}
 	Let $\A$ and $\B$ be a algebras of the same signature $\tau$. If $\B\in \HSPfin(\A)$, then the natural homomorphism from $\Clo(\A)$ onto $\Clo(\B)$ exists and is continuous. 
\end{prop}
 \begin{proof}
 	We show the statement for the cases where $\B$ is a finite product of $\A$, or a subalgebra of $\A$, or a homomorphic image of $\A$; the full statement then follows by combining the three. It is well-known that in all three cases, the natural homomorphism exists; this is because products, subalgebras, and homomorphic images of $\A$ satisfy at least the equations between $\tau$-terms that hold in $\A$. It thus remains to show that the natural homomorphism $\xi$ from $\Clo(\A)$ onto $\Clo(\B)$ is continuous.
	
	 Assume first that $\B=\A^n$ for some finite $n\geq 1$. Let $U$ be an open set from the subbasis of the topology on $\Clo(\A^n)\uk$, where $k\geq 1$; that is, there exist a $k$-tuple $a\in(A^n)^k$ and a value $v\in A^n$ such that $U$ consists precisely of those $k$-ary terms of $\A^n$ which send $a$ to $v$. Now viewing $a$ as a matrix in $A^{k\mult n}$, denote for all $1\leq i\leq n$ by $c_i$ the $i$-th column of $a$. Then $\xi^{-1}[U]$ consists of those $k$-ary terms of $\A$ which send $c_i$ to $v_i$, for all $1\leq i\leq n$; an open set. Hence, $\xi$ is continuous.
	 
	 Assume now that $\B$ is a subalgebra of $\A$. Then the preimage of any subbasis set $U$ of $\clob\uk$ is equal to $U$, and hence also open in $\cloa\uk$.
	 
	 Finally, let $\B$ be a homomorphic image of $\A$. Then $\B$ is isomorphic to $\A/{\sim}$ for a congruence relation ${\sim}$ of $\A$, and we may assume $\B=\A/{\sim}$. Let $U$ be a subbasis set of the topology of $\clob\uk$; so $U$ consists of those functions in $\clob\uk$ which send a certain tuple $a$ of $(\A/{\sim})^k$ to some $v\in (\A/{\sim})$. Then a function $f\in \cloa\uk$ is an element of $\xi\inv[U]$ iff there exists a $k$-tuple $c\in\A^k$ and $d\in \A$ such that $c_i\in a_i$ for all $1\leq i\leq n$, such that $d\in v$, and such that $f(c)=d$. For fixed $c,d$, the set of all $f\in \cloa\uk$ satisfying $f(c)=d$ is an open set in $\cloa\uk$, and so $\xi\inv[U]$ is the union of open sets and itself open.
	\end{proof}

\subsection{The converse}

We will now show that (1) implies (2) in Theorem~\ref{thm:topo-birk}.

Let $X,Y$ be countably infinite sets, and let $G$ be a group acting on $Y$. We equip the set $Y^X$ of all functions from $X$ to $Y$ with the topology of the Baire space, i.e., we consider $Y$ as a discrete space and give $Y^X$ the product topology.
Now define an equivalence relation ${\sim}_G$ on $Y^X$ which identifies two functions $f,g\in Y^X$ iff there exists $\alpha\in G$ such that $f=\alpha g$. We then consider the factor space $Y^X/{\sim}_G$ with the quotient topology, and write also $Y^X/G$ for this space; therefore, a subset $O\subseteq Y^X/ G$ is open iff $\bigcup O$ is open in $Y^X$.

\begin{prop}\label{prop:compact}
	Let $X,Y$ be countably infinite sets, and let $G$ be a group which acts on $Y$. Then $Y^X/G$ is compact if and only if the action of $G$ on $Y$ is oligomorphic.
\end{prop}
\begin{proof}
	We first prove that if the action of $G$ is oligomorphic, then $Y^X/G$ is compact. Say without loss of generality $X=\omega$. Pick for every $n\geq1$ and every orbit of the componentwise action of $G$ on $Y^n$ a representative tuple of this orbit in such a way that being a representative of an orbit is closed under taking initial segments; this can be done inductively. Write $R$ for the set of representatives. When we partially order $R$ by saying for $a,b\in R$ that $a$ is smaller or equal than $b$ if and only if $a$ is an initial segment of $b$, then $R$ becomes a finitely branching tree, the branches of which are elements of $Y^X$.
	 Consider the subspace $B$ of $Y^X$ of those functions which are branches of $R$; in other words, for $f\in Y^X$ we have $f\in B$ if and only if the restriction $f\rest_n$ of $f$ to $\{0,\ldots,n\}$ is in $R$, for all $n\geq 1$. Then $B$ is compact by Tychonoff's theorem as it is homeomorphic with a closed subspace of $\prod_{n\in\omega} k(n)$, where $k(n)$ is the (finite) number of representatives of length $n$. Moreover, $G\cdot B:=\{\alpha f\,|\, \alpha\in G\wedge f\in B\}$ is dense in $Y^X$, and so $(G\cdot B) / G$ is dense in $Y^X/G$. But no two elements $f,g$ of $B$ satisfy $f \sim_G g$, and so $(G\cdot B) / G$ is homeomorphic to $B$. Hence, $Y^X/G$ has a dense compact subset, proving that $Y^X/G$ is compact itself.
	 
	 For the other direction, assume that the action of $G$ is not oligomorphic. Pick an $n\geq 1$ such that the componentwise action of $G$ on $Y^n$ has infinitely orbits, and enumerate these orbits by $(O_i)_{i\in \omega}$. Now for all $i\in \omega$, let $U_i$ consist of all classes $[f]_{{\sim}_G}$ in $Y^X/G$ with the property that $f\rest_n$ belongs to $O_i$; this is well-defined since for all $f,g\in Y^X$ with $f \sim_G g$ we have that $f\rest_n$ belongs to $O_i$ iff $g\rest_n$ belongs to $O_i$. Then $Y^X/G$ is the disjoint union of the $U_i$. But each $U_i$ is open, and hence $Y^X/G$ is not compact.
\end{proof}

We remark that the space $Y^X/G$ is not Hausdorff, which explains that it can have a dense compact subset which is not equal to the whole space -- some readers might have wondered about this.

\begin{lem}\label{lem:invariantSubspace}
	Let $X,Y$ be countable sets, and let $G$ be a group with an oligomorphic action on $Y$. Let $S$ be a closed subset of $Y^X$ which is invariant under $G$, i.e., $G\cdot S\subseteq S$. Then $S/G$ is compact.
\end{lem}
\begin{proof}
	$S/G$ is a closed subspace of the compact space $Y^X/G$.
\end{proof}

For a structure $\Delta$, we write $\Emb(\Delta)$ for the set of self-embeddings of $\Delta$.

\begin{cor}\label{cor:polCompact}
 	Let $\Delta$ be an $\omega$-categorical structure. Then the following spaces are compact:
	\begin{itemize}
	\item  $\Emb(\Delta)/\Aut(\Delta)$;
	\item $\End(\Delta)/\Aut(\Delta)$;
	\item $\pol^{(k)}(\Delta)/\Aut(\Delta)$, for all $k\geq 1$.
		\end{itemize}
		Moreover, if $\A$ is a locally oligomorphic algebra and $\G$ is the group of all invertible unary bijections in $\overline{\Clo(\A)}$,  then $\ccloa\uk/\G$ is compact, for all $k\geq 1$.
\end{cor}
\begin{proof}
$\Emb(\Delta)$, $\End(\Delta)$ and $\pol^{(k)}(\Delta)$ are closed subsets of $\Delta^\Delta$ and $\Delta^{\Delta^k}$, respectively, which are invariant under $\Aut(\Delta)$. Since $\Delta$ is $\omega$-categorical, the action of $\Aut(\Delta)$  on $\Delta$ is oligomorphic by the theorem of Engeler, Svenonius, and Ryll-Nardzewski (see~\cite{HodgesLong}), and hence the first statement follows Lemma~\ref{lem:invariantSubspace}. The argument for the second statement is identical.
\end{proof}

Note that $\pol(\Delta)/\Aut(\Delta)$ is never compact since it is the disjoint union of the spaces $\pol^{(k)}(\Delta)/\Aut(\Delta)$.

\begin{nota}
	Let $D$ be a set, and let $f$ be a $k$-ary function on $D$ for some $k\geq 1$. If $C\in D^{m\mult k}$ for some $m\geq 1$, then we write $f(C)$ for the tuple of size $m$ obtained by applying $f$ to each row of the matrix $C$.
\end{nota}

\begin{lem}\label{lem:matrix}
	Let $\A, \B$ algebras of the same signature, where $\A$ is locally oligomorphic. Assume that the natural homomorphism $\xi$ from $\Clo(\A)$ onto $\Clo(\B)$ exists and is continuous. Then for all finite $F\subseteq B$ and all $k\geq 1$ there exist an $m\geq 1$ and $C\in A^{m\mult k}$ such that for all $f,g\in\Clo^{(k)}(\A)$ we have that $f(C)= g(C)$ implies $\xi(f)\rest_F= \xi(g)\rest_F$.
	\end{lem}
	\begin{proof}
	We denote the unique continuous extension of $\xi$ to $\ccloa$ by $\bar \xi$. So $\bar \xi$ is a continuous mapping from $\ccloa$ into $\cclob$. Moreover, it is a homomorphism: if $f\in \ccloa^{(n)}$ and $g_1,\ldots,g_n\in\ccloa^{(l)}$, where $n,l\geq 1$, then there exist sequences $(f^i)_{i\in\omega}$ and $(g_j^i)_{i\in\omega}$ of functions in  $\cloa^{(n)}$ and in  $\cloa^{(l)}$ which converge to $f$ and $g_j$, respectively, and so 
	\begin{align*}
	\bar\xi(f(g_1,\ldots,g_n))&=\bar\xi(\lim_{i\To\infty}(f^i(g_1^i,\ldots,g_l^i)))=\lim_{i\To\infty}\xi(f^i(g_1^i,\ldots,g_n^i))\\&=
	\lim_{i\To\infty}\xi(f^i)(\xi(g_1^i),\ldots,\xi(g_n^i))=\bar\xi(f)(\bar\xi(g_1),\ldots,\bar\xi(g_n)).
	\end{align*}
	We will prove the existence of $m\geq 1$ and $C\in A^{m\mult k}$ such that for all $f,g\in\ccloa\uk$ we have that $f(C)= g(C)$ implies $\bar\xi(f)\rest_F= \bar\xi(g)\rest_F$; the lemma then clearly follows.

Recall that the basic open sets of $\ccloa\uk$ are precisely the sets of the form
		$$
			O_{D,a}:=\{f \in\ccloa\uk\,|\, f(D)= a\},
		$$
		for $l\geq 1$, a matrix $D\in A^{l\mult k}$ and a vector $a\in A^l$; the basic open sets of $\cclob\uk$ are defined similarly. Call a basic open set $O$ of $\ccloa\uk$ an \emph{island} iff $\bar\xi(f)\rest_F=\bar\xi(g)\rest_F$ for all $f,g\in O$. From the definition of the basic open sets it is clear that for $f\in\ccloa\uk$, the set of all $h\in \overline{\Clo(\B)}^{(k)}$ which agree with $\bar\xi(f)$ on $F$ is open in $\overline{\Clo(\B)}$. Hence, the continuity of $\bar\xi$ implies that every $f\in\ccloa\uk$ is contained in a basic open island.
		
				Write $\G$ for the group of unary invertible bijections in $\overline{\Clo(\A)}$. Then $\G$ is oligomorphic as $\A$ is locally oligomorphic. 
		Observe next that for any basic open island $O_{D,a}$ of $\ccloa\uk$, the set $\G\cdot O_{D,a}=\{\alpha f\,|\, \alpha\in\G\wedge f\in O_{D,a}\}$ is an open subset of $\ccloa\uk$ which is invariant under $\G$; hence, it defines an open subset $V_{D,a}$ of $\ccloa\uk/\G$, namely the set of all ${\sim}_\G$-classes which have a representative in $O_{D,a}$. So every class $[f]_{{\sim}_\G}$ is contained in some set $V_{D,a}$ for a basic open island $O_{D,a}$. 	
		Since $\ccloa\uk/\G$ is compact by Proposition~\ref{prop:compact}, there are finitely many basic open islands $O_{D_1,a_1},\ldots,O_{D_n,a_n}$ such that the corresponding sets $V_{D_i,a_i}$ cover $\ccloa\uk/\G$. We then have that $\ccloa\uk$ is covered by the sets $\G\cdot O_{D_i,a_i}$. Set $m:=l_1+\cdots+l_n$, where $l_i$ denotes the number of rows of $D_i$, for $1\leq i\leq n$. Let $C$ be the matrix of dimension $m\mult k$ which is obtained by superposing the $D_i$. To see that $C$ satisfies the desired property, let $f,g\in\ccloa\uk$. Assume wlog that $f\in \G\cdot O_{D_1,a_1}$; then there exists $\alpha\in\G$ such that $f(D_1)=\alpha(a_1)$. Since $f(C)=g(C)$, we have $f(D_1)=g(D_1)$, and so also $g(D_1)=\alpha(a_1)$. Hence, $\alpha^{-1}f$ and $\alpha^{-1}g$ are in $O_{D_1,a_1}$, implying $\bar\xi(\alpha^{-1}f)\rest_F=\bar\xi(\alpha^{-1}g)\rest_F$ since $O_{D_1,a_1}$ is an island. Thus, $\bar\xi(f)\rest_F=\bar\xi(g)\rest_F$ since $\bar\xi$ is a homomorphism.
	\end{proof}

\begin{definition}
	We say that an algebra $\A$ of signature $\tau$ is \emph{finitely generated} iff there exists a finite subset $F$ of the domain of $\A$ such that the only subalgebra  of $\A$ containing $F$ is $\A$ itself; in other words, every element $a$ of $\A$ can be written as $t^\A(b_1,\ldots,b_k)$ for some $k\geq 1$, a $k$-ary $\tau$-term $t$, and $b_1,\ldots,b_k\in F$.
\end{definition}

\begin{prop}\label{prop:converse}
	Let $\A, \B$ be algebras of the same signature $\tau$, where $\A$ is locally oligomorphic and $\B$ is finitely generated. If the natural homomorphism from $\Clo(\A)$ onto $\Clo(\B)$ exists and is continuous, then $\B\in \HSPfin(\A)$.
\end{prop}
\begin{proof}
	 Let $F=\{b_1,\ldots,b_k\}$ be a set of generators of $\B$, and let $m\geq 1$ and $C\in A^{m\mult k}$ be given by Lemma~\ref{lem:matrix}. Let $\SSSS$ be the subalgebra of $\A^m$ generated by the columns $c_1,\ldots,c_k$ of $C$; so the elements of $\SSSS$ are precisely those of the form $t^{\A^m}(c_1,\dots,c_k)$, for a $k$-ary $\tau$-term $t$. Define a function $\mu\colon \SSSS\To \B$ by setting $\mu(t^{\A^m}(c_1,\dots,c_k)):=t^\B(b_1,\ldots,b_k)$. Then $\mu$ is well-defined, for if $t^{\A^m}(c_1,\dots,c_k)=s^{\A^m}(c_1,\dots,c_k)$, then $t^\B(b_1,\ldots,b_k)=s^\B(b_1,\ldots,b_k)$ by the properties of $C$. Since $\B$ is generated by $F$, $\mu$ is onto. We claim that $\mu$ is moreover a homomorphism; it then follows that $\B$ is the  homomorphic image of the subalgebra $\SSSS$ of $\A^m$, and so $\B\in\HSPfin(\A)$. To this end, let $f$ be any function symbol of $\tau$, let $n$ be its arity, and let $s_1,\ldots,s_n\in\SSSS$. Write $s_i=t_i^{\A^m}(c_1,\ldots,c_k)=t_i^\SSSS(c_1,\ldots,c_k)$ for all $1\leq i\leq n$. Then
	 \begin{align*}
	 	\mu(f^\SSSS(s_1,\ldots,s_n))&=\mu(f^\SSSS(t_1^\SSSS(c_1,\ldots,c_k),\ldots,t_n^\SSSS(c_1,\ldots,c_k)))=
		\mu(f^\SSSS(t_1^\SSSS,\ldots,t_n^\SSSS)(c_1,\ldots,c_k))\\
		&=\mu((f(t_1,\ldots,t_n))^\SSSS(c_1,\ldots,c_k))=
		(f(t_1,\ldots,t_n))^\B(b_1,\ldots,b_k)\\
		&=
		f^\B(t_1^\B(b_1,\ldots,b_k),\ldots,t_n^\B(b_1,\ldots,b_k))\\&=
		f^\B(\mu(s_1),\ldots,\mu(s_n)).
	 \end{align*}
	 
	 \end{proof}

\begin{prop}\label{prop:finitelygenerated}
	Let $\B$ be an algebra which is locally oligomorphic. Then $\B$ is finitely generated.
\end{prop}
\begin{proof}
	Let $\G$ be the permutation group of invertible unary bijections of $\cclob$. Since $\B$ is locally oligomorphic, the action of $\G$ on $B$ has finitely many orbits. Picking a representative from each orbit one obtains a generating set for $\B$.
\end{proof}
 Theorem~\ref{thm:topo-birk} now follows from Propositions~\ref{prop:hin}, \ref{prop:converse}, and \ref{prop:finitelygenerated}. Note that in the theorem, it would have been sufficient to assume that $\B$ be finitely generated rather than locally oligomorphic, but since we are mainly interested in polymorphism clones of $\omega$-categorical structures we have chosen to formulate the theorem as it is. The following is the stronger variant which follows from Propositions~\ref{prop:hin} and \ref{prop:converse}.
 
  \begin{theorem}\label{thm:variant}
 Let $\A,\B$ be algebras with the same signature, where $\A$ is locally oligomorphic and $\B$ is finitely generated. 
Then the following three statements are equivalent. 
\begin{enumerate}
\item The natural homomorphism from $\Clo(\A)$ onto $\Clo(\B)$ exists and is continuous. 
\item $\B \in \HSPfin(\A)$.
\item $\B$ is contained in the pseudovariety generated by $\A$.
\end{enumerate}
\end{theorem}

\section{Pseudovariety Examples}
\label{sect:examples1}

We now give two examples examining the continuity condition on the natural homomorphism in Theorem~\ref{thm:topo-birk}.
The first example is due to Keith Kearnes~\cite{KeithExample}, and demonstrates that there are oligomorphic algebras $\A$ such that the variety generated by $\A$ contains finite members which the pseudovariety generated by $\A$ does not contain.
 
\begin{prop}\label{prop:ex1:1}
There are algebras $\A, \B$ with common signature such that 
\begin{itemize}
\item $\A$ is locally oligomorphic;
\item $\B$ is finite;
\item $\B\in\HSP(\A)$;
\item $\B\nin\HSPfin(\A)$.
\end{itemize} 
Hence, the natural homomorphism from $\Clo(\A)$ onto $\Clo(\B)$ exists but is not continuous.
\end{prop}
\begin{proof}
	Let the signature $\tau$ consist of unary function symbols $(f_i)_{i\in\omega}$ and $(g_i)_{i\in\omega}$. Let $\A$ be any algebra on $\omega$ with signature $\tau$ such that the functions ${f_i^\A}$ form a locally oligomorphic permutation group, such that no $g_i^\A$ is injective, and such that $f_0^\A$ is contained in the topological closure of $\{g_i^\A\}_{i\in\omega}$. Let $\B$ be the $\tau$-algebra on $\{0,1\}$ such that $f_i^\B$ is the identity function for all $i\in\omega$ and such that $g_i^\B$ is the constant function with value $0$. It is easy to see that the natural homomorphism from $\cloa$ onto $\clob$ exists. However, it is not continuous since $f_0^\A$ is contained in the topological closure of $\{g_i^\A\}_{i\in\omega}$, but $f_0^\B$ is not contained in the topological closure of $\{g_i^\B\}_{i\in\omega}$.
\end{proof}

We remark that one can easily modify the previous example to obtain algebras $\A, \B$ with finite signature and the same properties. On the other hand, by taking an uncountable signature, one can make $\A$ even oligomorphic.

The next example becomes relevant when one has (concrete) clones without a signature of a corresponding algebra; this is for example the case for polymorphism clones of structures, as in Theorem~\ref{thm:main} and in the following sections. It shows that when we are given two such clones $\CC,\D$, then it might happen that there exists a homomorphism from $\CC$ onto $\D$ which is not continuous, as well as a continuous clone homomorphism onto $\D$. In other words, when we make algebras out of $\CC$ and $\D$ by matching the functions in $\CC$ and $\D$ with an appropriate functional signature $\tau$, then we might do so in such a way that the natural homomorphism from $\CC$ onto $\D$ exists and is continuous, and in another way such that  the natural homomorphism from $\CC$ onto $\D$ exists but is not continuous.

\begin{prop}
There are algebras $\A, \B$ such that 
\begin{itemize}
\item $\A$ is oligomorphic;
\item $\B$ is finite;
\item there exists a non-continuous clone homomorphism from $\cloa$ onto $\clob$;
\item there exists a continuous clone homomorphism from $\cloa$ onto $\clob$.
\end{itemize}
	\end{prop}
\begin{proof}
	Let $\B$ be as in Proposition~\ref{prop:ex1:1}. Let $\A$ be the algebra on $\omega$ which has the following three sets of functions.
	\begin{align*}
		\F_1&:=\{f\in\omega^\omega \,|\, f(0)=f(1)=1\text{ and } (\forall n\geq 2\; f(n)\geq 2)\text{ and } f \text{ is not surjective}\},\\
		\F_2&:=\{f\in\omega^\omega \,|\, f(0)=f(1)=1\text{ and }  f{\rest}_{[2,\infty)} \text{ is a permutation on }  [2,\infty]\},\\
		\F_3&:=\{f\in\omega^\omega \,|\, f(0)=0\text{ and } f(1)=1\text{ and } f \text{ is a permutation on }\omega\}.
\end{align*}

	Now observe that if $f\in \F_i$ and $g\in\F_j$, then $f\circ g\in \F_{\min(i,j)}$. The function which sends all elements of $\F_1\cup \F_2$ to the constant function of $\B$ and all elements of $\F_3$ to the identity induces a continuous homomorphism from $\Clo(\A)$ onto $\Clo(\B)$. On the other hand, the function which sends all elements of $\F_1$ to the constant unary function of $\clob$ and all elements of $\F_2\cup\F_3$ to the identity in $\clob$ induces a non-continuous homomorphism from $\cloa$ onto $\clob$.
\end{proof}
 
\section{Primitive Positive Interpretations}
\label{sect:interpretations}
In this section we prove Theorem~\ref{thm:main}. 
Our definition of interpretations follows~\cite{HodgesLong} and is standard, and will be recalled in the following.
Let $\tau$ be a signature, and let $\Gamma$ be a $\tau$-structure.
If $\delta(x_1,\dots,x_k)$ is a first-order $\tau$-formula with $k$ free variables $x_1,\dots,x_k$, we write $\delta(\Gamma^k)$ for the $k$-ary relation that is defined by $\delta$ on $\Gamma$.

An atomic $\tau$-formula is called \emph{unnested} iff it is of the form $x_0=x_1$, of the form $x_0 = f(x_1,\dots,x_n)$, or of the form $R(x_1,\dots,x_n)$, for some $n$-ary function symbol $f \in \tau$ or relation symbol $R \in \tau$, and variables $x_0,x_1,\dots,x_n$.
It is straightforward to see that every atomic $\tau$-formula is
equivalent to a primitive positive $\tau$-formula whose atomic
subformulas are unnested (see Theorem~2.6.1 in~\cite{HodgesLong}). 

\begin{definition}
A $\sigma$-structure $\Delta$ has a \emph{(first-order) interpretation $I$} in a $\tau$-structure $\Gamma$ iff there exists a natural number $d \geq 1$, called the \emph{dimension} of $I$, and
\begin{itemize}
\item a $\tau$-formula $\delta_I(x_1, \dots, x_d)$ -- called \emph{domain formula},
\item for each unnested atomic $\sigma$-formula $\phi(y_1,\dots,y_k)$ a $\tau$-formula $\phi_I(\overline x_1, \dots, \overline x_k)$ where the $\overline x_i$ denote disjoint $d$-tuples of distinct variables -- called the \emph{defining formulas},
\item a surjective map $h \colon \delta_I(\Gamma^d) \rightarrow \Delta$ -- called \emph{coordinate map},
\end{itemize}
such that for every unnested atomic $\sigma$-formula $\phi$ and all tuples $\overline a_i \in \delta_I(\Gamma^d)$ 
\begin{align*}
\Delta \models \phi(h(\overline a_1), \dots, h(\overline a_k)) \; 
& \Leftrightarrow \; 
\Gamma \models \phi_I(\overline a_1, \dots, \overline a_k) \; .
\end{align*}
\end{definition}
If the formulas $\delta_I$ and $\phi_I$ are primitive positive (existential positive), 
we say that the interpretation $I$ is \emph{primitive positive} (\emph{existential positive}). Note that the dimension $d$, the set $S := \delta_I(\Gamma^d)$, and the coordinate map $h$ determine the defining formulas up to logical equivalence; hence, we sometimes denote an interpretation by $I = (d,S,h)$.

\subsection{Primitive positive interpretations and pseudovarieties} For $\omega$-categorical structures $\Gamma$, 
primitive positive interpretability in $\Gamma$ can be characterized in terms of the pseudovariety generated by a polymorphism algebra of $\Gamma$. 
Via the results of the previous section, pseudovarieties also correspond to topological
clones -- so they provide the link between primitive positive interpretations 
and topological clones, which will be used to prove Theorem~\ref{thm:main} 
in Section~\ref{ssect:interpret-topoclones} (confer also Figure~\ref{fig:diagram}).

\begin{definition}
Let $\Gamma$ be a structure, and $\A$ an algebra. 
Then $\A$ is called a \emph{polymorphism algebra} of $\Gamma$ iff $\A$ and $\Gamma$ have the same domain, and the set of operations of $\A$ is precisely 
the set of polymorphisms of $\Gamma$. 
\end{definition}
Clearly, every structure $\Gamma$ has a polymorphism algebra, which can be obtained by assigning function names to the polymorphisms in some arbitrary way.

\begin{thm}\label{thm:pp-pseudovar} 
  Let $\Gamma$ be a finite or $\omega$-categorical structure, and let $\Delta$ be an arbitrary structure. 
  Then the following are equivalent.
  \begin{enumerate}
 	\item for every polymorphism algebra $\C$ of $\Gamma$ there is an algebra $\B \in \HSPfin(\C)$ such that $\Clo(\B) \subseteq \Pol(\Delta)$;
  \item there is a polymorphism algebra $\C$ of $\Gamma$
  and an algebra $\B \in \HSPfin(\C)$ such that $\Clo(\B)\subseteq \Pol(\Delta)$;
 	\item  $\Delta$ has a primitive positive interpretation
  in $\Gamma$.
  \end{enumerate}
\end{thm}

The equivalence between (1) and (2) emphasizes the fact that for our purposes, it does not matter in what way we assign function names to the polymorphisms of $\Gamma$. Theorem~\ref{thm:pp-pseudovar} already appeared in the survey 
article~\cite{BodirskySurvey}; it has been inspired
by results obtained in the context of constraint satisfaction problems for finite structures~\cite{JBK}.
Since we need Theorem~\ref{thm:pp-pseudovar} in a more detailed form 
(Proposition~\ref{prop:pp-pseudovar-detail}),
we provide its full proof here. 

Let $\Gamma$ be $\tau$-structure with domain $D$, and $R \subseteq D^k$ a $k$-ary relation. We say that $R$ is \emph{primitive positive definable} in $\Gamma$ iff there exists a primitive positive $\tau$-formula $\phi(x_1,\dots,x_k)$ such that for all $(c_1,\dots,c_k) \in D^k$
it is true that $(c_1,\dots,c_k) \in R$ if and only if $\Gamma$ satisfies $\phi(c_1,\dots,c_k)$. We say that a $\tau$-formula $\phi$ with $k$ free variables is \emph{preserved} 
by a function $f \colon D^l \rightarrow D$ (over $\Gamma)$ iff 
for all $t_1^1,\dots,t_l^k \in D$, if $\Gamma \models 
\phi(t^1_i,\dots,t^k_i)$ for all $i \leq l$, then $\Gamma \models
\phi(f(t_1^1,\dots,t_l^1),\dots,f(t^k_1,\dots,t^k_l))$. 
Note that $f$ is a polymorphism of $\Gamma$ if and only if
$f$ preserves all atomic unnested $\tau$-formulas over $\Gamma$. 
We say that a relation $R\subseteq D^k$ (a function $g\colon D^k\To D$) is preserved by $f$ iff $f$ is a polymorphism of the structure $(D;R)$ (of $(D;g)$).

We need the following characterization of primitive positive definability in $\omega$-categorical structures $\Gamma$; for finite structures $\Gamma$, this is due to~\cite{Geiger,BoKaKoRo}. 

\begin{theorem}[from~\cite{BodirskyNesetrilJLC}]\label{thm:pp-def}
Let $\Gamma$ be finite or $\omega$-categorical. Then a relation $R$ has a primitive positive definition in $\Gamma$ if and only if $R$ is preserved by all polymorphisms of $\Gamma$.
\end{theorem}

For example, when $D$ is the domain of an $\omega$-categorical
structure $\Gamma$ and $\C$ is a polymorphism algebra of $\Gamma$, then an equivalence relation relation $R \subseteq D^2$ is a congruence of $\C$ if and only if $R$ is primitive positive definable in $\Gamma$. 

\begin{proof}[Proof of Theorem~\ref{thm:pp-pseudovar}]
The implication from $(1)$ to $(2)$ follows from the existence
of a polymorphism algebra $\C$ of $\Gamma$. 

$(2) \Rightarrow (3)$. 
Write $\tau$ for the signature of $\C$.
 There exists a finite number $d \geq 1$,
 a subalgebra $\bf S$ of $\C^d$ with domain $S$,
 and a surjective homomorphism $h$ 
from $\bf S$ to $\bf B$.
We claim that $\Delta$ has the primitive positive interpretation $I:=(d,S,h)$ in $\Gamma$. All operations of $\C$ preserve $S$ (viewed as a $d$-ary relation 
over $\Gamma$), since ${\bf S}$ is a subalgebra of $\C^d$. Theorem~\ref{thm:pp-def} implies that $S$ has a primitive positive definition $\delta(x_1,\dots,x_d)$ in $\Gamma$, which becomes the 
domain formula $\delta_I$. 

Let $\psi$ be an unnested atomic formula over the signature of $\Delta$ and with $k$ free variables $x_1,\dots,x_k$. 
Let $R \subseteq C^{dk}$ be the relation defined by
$$(a^1_1,\dots,a^d_1,\dots, a^1_k,\dots,a^d_k) \in R \; \Leftrightarrow 
\; \Delta \models \psi(h(a^1_1,\dots,a^d_1),\dots,h(a^1_k,\dots,a^d_k)) \; ,$$
and let $f \in \tau$ be arbitrary. 
By assumption, $f^\B$ preserves $\psi$. 
Since $h$ is a homomorphism, it follows that $f^\C$ 
preserves $R$. 
We conclude that all polymorphisms of $\Gamma$
preserve $R$. Since $\Gamma$ is $\omega$-categorical and by Theorem~\ref{thm:pp-def}, the
relation $R$ has a primitive positive definition in $\Gamma$,
which becomes the defining formula for $\psi(x_1,\dots,x_k)$.
So $I$ is indeed a
primitive positive interpretation of $\Delta$ in $\Gamma$.

To prove $(3) \Rightarrow (1)$, 
suppose that $\Delta$ has a primitive positive interpretation $I = (d,S,h)$
in $\Gamma$. Let $\C$ be a polymorphism algebra of $\Gamma$, and let $\tau$ be the signature of $\C$.
We have to show that $\HSPfin(\C)$ contains a $\tau$-algebra ${\bf B}$ such
that all operations in ${\bf B}$ are polymorphisms of $\Delta$.
The set $S$ is preserved by all operations of $\Clo(\C) = \Pol(\Gamma)$, because it is
primitive positive definable in $\Gamma$ by the domain formula of $I$ (Theorem~\ref{thm:pp-def}). Therefore, $S$ induces a subalgebra $\bf S$ of $\C^d$.
Let $K$ be the kernel of the coordinate map $h$
of $I$. Then for all tuples $\overline a, \overline b \in S$,
the $2d$-tuple $(\overline a,\overline b)$
satisfies $=_I$ in $\Gamma$ if and only if
$(\overline a,\overline b) \in K$. Since $=_I$ is primitive positive
definable in $\Gamma$, it is preserved by all operations of $\C$ by
Theorem~\ref{thm:pp-def}.  It follows that $K$ is a congruence of
$\SSSS$. As a consequence, $h$ induces a $\tau$-algebra $\B$ on its image, which equals the domain of $\Delta$, in such a way that $h$ is a homomorphism from $\SSSS$ onto $\B$: let $f \in \tau$ be $m$-ary, and let $c_1,\dots,c_m$ be arbitrary elements of $\Delta$. Then pick $\overline a_1,\dots,\overline a_m \in S$ such that $h(\overline a_i) = c_i$, and define $f^{\B}(c_1,\dots,c_m) := h(f^\SSSS(\overline a_1),\dots,f^{\SSSS}(\overline a_m))$. This is well-defined since the kernel $K$ of $h$ is a congruence of $\SSSS$, and by definition of $\B$, $h$ is a homomorphism from $\SSSS$ onto $\B$. It remains to verify that for all $f\in\tau$, $f^{\bf B}$ is a polymorphism of $\Delta$,
i.e., every unnested atomic formula $\phi$ over $\Delta$ is preserved by $f^{\bf B}$. From the definitions of $\phi_I$ and $f^\B$, one easily sees that  $f^\B$ preserves $\phi$ over $\Delta$ if and only if
$f^{\bf C}$ preserves $\phi_I$ over $\Gamma$.
Since $f^{\bf C}$ is a polymorphism of $\Gamma$, 
and since $\phi_I$
is a primitive positive $\tau$-formula over $\Gamma$, $f^\C$ indeed
preserves $\phi_I$, and hence $f^\B$ preserves $\phi$. 
\end{proof}

The proof of Theorem~\ref{thm:pp-pseudovar} above gives
more information about the link between polymorphism algebras
and primitive positive interpretations, and we state it explicitly.  

\begin{prop}\label{prop:pp-pseudovar-detail}
  Let $\Gamma$ be a finite or $\omega$-categorical structure with domain $D$, and let $\Delta$ be an arbitrary structure with domain $B$. Then for all $d \geq 1$, $S \subseteq D^d$, and $h \colon S \rightarrow B$ the following are equivalent.
  \begin{enumerate}
  \item For every polymorphism algebra $\C$ of $\Gamma$ the set $S$ induces a subalgebra $\SSSS$ of $\C^d$, 
  the kernel of $h$ is a congruence of $\SSSS$, and
  the homomorphic image $\B$ of $\SSSS$ under $h$ satisfies $\Clo(\B) \subseteq \Pol(\Delta)$;
 	\item  $\Delta$ has the primitive positive interpretation $(d,S,h)$ in $\Gamma$.
  \end{enumerate}
\end{prop}

\subsection{Primitive positive interpretations and topological clones}\label{ssect:interpret-topoclones}
We can now show the first part of Theorem~\ref{thm:main}.

\begin{proposition}\label{prop:interpret-clone}
Let $\Gamma$ be finite or $\omega$-categorical, and $\Delta$ be arbitrary. Then $\Delta$ has a primitive positive interpretation in $\Gamma$ if and only if $\Delta$ is the reduct of a finite or $\omega$-categorical structure $\Delta'$ such that there exists a continuous clone homomorphism from $\Pol(\Gamma)$ to $\Pol(\Delta')$ whose image is dense in $\Pol(\Delta')$.
\end{proposition}
\begin{proof}
Let $\C$ be a polymorphism algebra of $\Gamma$.

Suppose first that $\Delta$
has a primitive positive interpretation in $\Gamma$.
By Theorem~\ref{thm:pp-pseudovar} there is an algebra $\B$
in the pseudovariety generated by $\C$ 
such that all operations of $\B$ are polymorphisms of $\Delta$.
Since $\Gamma$ is finite or $\omega$-categorical, $\C$ is finite or oligomorphic, and the algebra $\B$ is finite or oligomorphic as well. By Theorem~\ref{thm:topo-birk} the natural homomorphism $\xi$ from $\Clo(\C)$ 
onto $\Clo(\B)$ exists and continuous. 
Let $\Delta'$ be the structure
with the same domain as
$\B$ that contains all relations and all functions preserved by all
operations of $\B$. Since $\Clo(\B) \subseteq \Pol(\Delta')$, it
follows that $\Delta'$
is finite or $\omega$-categorical by the theorem of Engeler,
Svenonius, and Ryll-Nardzewski. Moreover, it is easy to see and
well-known that $\Pol(\Delta') = \overline{\Clo(\B)}$, so the image of
$\xi$ is dense in $\Pol(\Delta')$.
Since all operations of $\B$ are polymorphisms of $\Delta$,
all relations and functions of $\Delta$ are relations and functions of $\Delta'$, and
this shows that $\Delta$ is indeed a reduct of $\Delta'$.  

To prove the converse, let $\Delta'$ be a finite or $\omega$-categorical
structure such that $\Delta$ is a reduct of $\Delta'$, and such
that there is a continuous homomorphism $\xi$ from $\Pol(\Gamma)$ to $\Pol(\Delta')$ whose image is dense in $\Pol(\Delta')$. Let $\B$ be the algebra with the same domain as $\Delta$, the same signature $\tau$ as $\C$, and where $f \in \tau$ denotes the operation $\xi(f^\C)$ of $\Pol(\Delta')$.  Then $\overline{\Clo(\B)}=\Pol(\Delta')$ since the image of $\xi$ is dense in $\Pol(\Delta')$. Hence, $\B$ is finite or locally oligomorphic since $\Delta'$ is finite or $\omega$-categorical. We can therefore apply Theorem~\ref{thm:topo-birk} to infer $\B \in \HSPfin(\C)$. 
By Theorem~\ref{thm:pp-pseudovar}, $\Delta'$ has a primitive positive
interpretation in $\Gamma$. It follows that in particular $\Delta$
has a primitive positive interpretation in $\Gamma$.
\end{proof}

In Section~\ref{sect:examples2} we will present an example showing that in Proposition 22 we cannot simply require the continuous clone homomorphism $\xi$ to be surjective. In particular, the image of a closed oligomorphic clone under a continuous homomorphism need not be closed.

How do we recognize whether
two structures $\Gamma$ and $\Delta$
have isomorphic topological polymorphism clones?

\begin{definition}
Two structures $\Gamma$ and $\Delta$ such that 
$\Gamma$ has a primitive positive interpretation in $\Delta$
and $\Delta$ has a primitive positive interpretation in $\Gamma$
are called \emph{mutually primitive positive interpretable}. 
\end{definition}

We will see in Section~\ref{sect:examples2}
that there are $\omega$-categorical structures $\Gamma$ and $\Delta$ that are mutually primitive positive interpretable
and have non-isomorphic topological polymorphism clones. To characterize
the situation where $\Gamma$ and $\Delta$ have isomorphic topological polymorphism clones, we need the following stronger notion. 

\begin{definition}
Two structures $\Gamma$ and $\Delta$ are
called \emph{primitive positive bi-interpretable}\footnote{Here we follow the analogous definition for \emph{first-order bi-interpretability} as introduced in~\cite{AhlbrandtZiegler}.} 
iff there is an interpretation $I = (d_1,S_1,h_1)$ of $\Delta$ in $\Gamma$ and an interpretation $J = (d_2,S_2,h_2)$ of $\Gamma$ in $\Delta$
such that the $(1 + d_1 d_2)$-ary relation $R_{IJ}$ defined by
$$x = h_1(h_2(y_{1,1},\dots,y_{1,d_2}),\dots,h_2(y_{d_1,1},\dots,y_{d_1,d_2})) $$
is primitive positive definable in $\Delta$,
and the $(1 + d_1 d_2)$-ary relation $R_{JI}$ defined by
$$x = h_2(h_1(y_{1,1},\dots,y_{1,d_1}),\dots,h_1(y_{d_2,1},\dots,y_{d_2,d_1})) $$ 
is primitive positive definable in $\Gamma$. 
\end{definition}
In the following, we write $h_1 \circ h_2$ for the function
defined by $$(y_{1,1},\dots,y_{1,d_2},\dots,y_{d_1,1},\dots,y_{d_1,d_2}) \mapsto h_1(h_2(y_{1,1},\dots,y_{1,d_2}),\dots,h_2(y_{d_1,1},\dots,y_{d_1,d_2}))\; .$$ 

\begin{proposition}\label{prop:bi-interpret}
Let $\Gamma$ and $\Delta$ be finite or $\omega$-categorical. Then the following are equivalent.
\begin{enumerate}
\item $\Pol(\Gamma)$ and $\Pol(\Delta)$ are isomorphic as topological clones. 
\item $\Gamma$ has a polymorphism algebra $\A$,
and $\Delta$ has a polymorphism algebra $\B$ such
that $\HSPfin(\A)=\HSPfin(\B)$.
\item $\Gamma$ and $\Delta$ are primitive positive bi-interpretable. 
\end{enumerate}
\end{proposition}
\begin{proof}
We prove $(1) \Rightarrow (2) \Rightarrow (3) \Rightarrow (1)$. 
Let $\A$ be a polymorphism algebra of $\Gamma$ 
with signature $\tau$, 
and suppose that
$\Pol(\Gamma)$ and $\Pol(\Delta)$ are isomorphic via a homeomorphism
$\xi$. Let $\B$ be the algebra with the same domain as $\Delta$ and signature $\tau$ such that $f^{\B} = \xi(f^{\A})$ for all $f \in \tau$. Then $\B$ is a polymorphism algebra of $\Delta$, 
and it follows from Theorem~\ref{thm:topo-birk}
that $\HSPfin(\A) = \HSPfin(\B)$. Thus $(1)$ indeed implies $(2)$. 

$(2) \Rightarrow (3)$. 
Suppose that $\Gamma$ has a polymorphism algebra $\A$
and $\Delta$ has a polymorphism algebra $\B$ such that
$\HSPfin(\A) = \HSPfin(\B)$. So there is a $d_1 \geq 1$,
a subalgebra ${\bf S}_1$ of $\A^{d_1}$, and a surjective
homomorphism $h_1$ from $\SSSS_1$ to $\B$.
Similarly, there is a $d_2 \geq 1$, 
a subalgebra $\SSSS_2$ of $\B^{d_2}$, and a surjective homomorphisms $h_2$ from ${\bf S}_2$ to $\A$.
By Proposition~\ref{prop:pp-pseudovar-detail}, $I := (d_1,S_1,h_1)$
is an interpretation of $\Delta$ in $\Gamma$,
and $J := (d_2,S_2,h_2)$ is an interpretation of $\Gamma$ in $\Delta$. Because the statement is symmetric it
suffices to show that $R_{IJ}$ is primitive positive definable
in $\Delta$. 
Theorem~\ref{thm:pp-def} asserts that this is equivalent to showing that $h_1 \circ h_2$ is preserved by all operations $f^\B$ of $\B$. 
So let  $k$ be the arity of $f^\B$, let $D$ be the domain of $\Delta$, and let $b^i = (b^i_{1,1},\dots,b^i_{d_1,d_2})$ be elements of
 $D^{d_1d_2}$, for $1\leq i \leq k$. Then 
 \begin{align*}
 f^{\B}((h_1 \circ h_2)(b^1),\dots,(h_1 \circ h_2)(b^k))\ =\ & 
 h_1\big (f^\A(h_2(b^1_{1,1},\dots,b^1_{1,d_2}),\dots,h_2(b^k_{1,1},\dots,b^k_{1,d_2})),\dots,\\
 & \quad \;  f^\A(h_2(b^1_{d_1,1},\dots,b^1_{d_1,d_2}),\dots,h_2(b^k_{d_1,1},\dots,b^k_{d_1,d_2}))\big ) \\
\ =\ & (h_1 \circ h_2)(f^\B(b^1,\dots,b^k)) \; .
\end{align*}

$(3) \Rightarrow (1)$. 
Suppose that $\Gamma$ and $\Delta$
are primitive positive bi-interpretable via an interpretation
$I=(d_1,S_1,h_1)$ of $\Delta$ in $\Gamma$
and an interpretation $J=(d_2,S_2,h_2)$ of $\Gamma$ in $\Delta$.
Let $\A$ be a polymorphism algebra of $\Gamma$, and $\B$ be a polymorphism algebra of $\Delta$. 
Then by Proposition~\ref{prop:pp-pseudovar-detail} $S_1$
induces an
algebra $\SSSS_1$ in $\A^{d_1}$ and
$h_1$ is a surjective homomorphism from $\SSSS_1$ to an algebra ${\B'}$ satisfying
$\Clo(\B') \subseteq \Pol(\Delta)$. 
Similarly, $S_2$ induces in $\B^{d_2}$ an algebra
$\SSSS_2$ 
and $h_2$ is a homomorphism from $\SSSS_2$ onto an algebra $\A'$
such that $\Clo(\A') \subseteq \Pol(\Gamma)$. 
By Theorem~\ref{thm:topo-birk} the natural homomorphisms
$\xi_1$
from $\Clo(\A)$ onto $\Clo(\B')$ and $\xi_2$ from $\Clo(\B)$ onto $\Clo(\A')$ exist and are continuous. We will verify that $\xi_2 \xi_1$ is the identity on $\Clo(\A)$; the proof that $\xi_1 \xi_2$ on $\Clo(\B)$ is the identity is analogous. It then follows that $\xi_1$ and $\xi_2$ are isomorphisms and homeomorphisms between $\Clo(\A)$ and $\Clo(\B)$. 

Write $\tau$ for the signature of $\A$. Write $\C$ for the $\tau$-algebra on the domain of $\A$ obtained by setting $f^\C:=(\xi_2\xi_1)(f^\A)$ for all $f\in\tau$. Let $f \in \tau$ be $k$-ary; we show $f^\C=f^\A$. Let $a_1,\dots,a_k \in \Gamma$ be arbitrary. Since $h_2 \circ h_1$ is surjective onto $\Gamma$, there are 
$b^i = (b^i_{1,1},\dots,b^i_{d_1,d_2}) \in \Gamma^{d_1d_2}$ 
such that $a_i = h_2 \circ h_1(b^i)$.
Then
\begin{align*}
f^{\C}(a_1,\dots,a_k) \; = \; & f^{\C}(h_2 \circ h_1(b^1),\dots,h_2 \circ h_1(b^k)) \\
= \; &h_2\big (f^{\B'}(h_1(b^1_{1,1},\dots,b^1_{d_1,1}),\dots,h_1(b^k_{1,1},\dots,b^k_{d_1,1})),\dots,\\
 & \quad \;  f^{\B'}(h_1(b^1_{1,d_2},\dots,b^1_{d_1,d_2}),\dots,h_1(b^k_{1,d_2},\dots,b^k_{d_1,d_2}))\big )\\
= \; & h_2 \circ h_1(f^\A(b^1,\dots,b^k)) \\
= \; & f^\A(h_2 \circ h_1(b^1),\dots,h_2 \circ h_1(b^k)) \\
 = \; & f^\A(a_1,\dots,a_k) 
\end{align*}
where the second and third equations hold since $h_2$ and $h_1$ are 
algebra homomorphisms, and the fourth equation holds
because $f^\A$ preserves $h_2 \circ h_1$: this follows
from Theorem~\ref{thm:pp-def} and the assumption that $R_{JI}$
is primitive positive definable in $\Gamma$. Hence, $f^\A=f^{\C}=\xi_2 \xi_1(f^\A)$ for all $f\in\tau$, which is what we had to show.
\end{proof}

The following fact has been proven recently for \emph{finite} algebras, independently by Markovi\'c, Maroti, McKenzie~\cite{MarkovicMarotiMcKenzie}, and by Davey, Jackson, Pitkethly, and Szab\'o~\cite{DaveyJacksonPitkethlySzabo}. 
An algebra ${\bf A}$ is called \emph{finitely related}
iff there exists a structure $\Gamma$ with the same domain as ${\bf A}$ 
and with finite relational signature such that $\overline{\Clo({\bf A})} = \Pol(\Gamma)$. We present a generalization to all locally oligomorphic algebras. 

\begin{corollary}
Let ${\bf A}$ and ${\bf B}$ be finite or locally oligomorphic algebras 
such that $\overline{\Clo(\A)}$ and $\overline{\Clo(\B)}$ are isomorphic as topological clones. 
Then ${\bf A}$ is finitely related if and only if ${\bf B}$ is finitely related. 
\end{corollary}
\begin{proof}
Suppose that $\A$ is finitely related; that is,
there exists a structure $\Gamma$ with finite relational signature
such that $\overline{\Clo({\A})} = \Pol(\Gamma)$.
Let $\Delta$ be the relational structure with the same domain as $\B$
that has all relations that are preserved by
all operations of $\B$. Then $\Pol(\Delta)=\overline{\Clo(\B)}$, and
thus it suffices to show that $\Delta$ has a reduct $\Delta'$ with
finite 
signature and the same polymorphisms as $\Delta$.

Note that the automorphisms of $\Gamma$ and $\Delta$ are exactly the
unary invertible operations in $\overline{\Clo(\A)}$ and
$\overline{\Clo(\B)}$, respectively.
Since $\A$ and $\B$ are finite or locally oligomorphic,
$\Gamma$ and $\Delta$ are finite or $\omega$-categorical.
By Proposition~\ref{prop:bi-interpret}, $\Gamma$ and
$\Delta$ are primitive positive bi-interpretable. Let $I_1$ and $I_2$
be the corresponding interpretations of $\Gamma$ in $\Delta$ and
$\Delta$ in $\Gamma$, respectively. Let $\sigma$ be the signature of
$\Delta$, and let $\sigma' \subseteq \sigma$ be the set of all relation symbols that appear in all the formulas of $I_1$; since the signature $\tau$
of $\Gamma$ is finite, $\sigma'$ is finite as well. Let $\Delta'$ be
the
$\sigma'$-reduct of $\Delta$.
We will show that there is a primitive positive definition of $\Delta$
in $\Delta'$; by  Theorem~\ref{thm:pp-def}, this implies that $\Delta$
and $\Delta'$ have the same polymorphisms.

Let $\psi$ be an atomic $\sigma$-formula with $k$ free variables
$x_1,\dots,x_k$. We specify an
equivalent primitive positive $\sigma'$-formula.
Suppose that the interpretation
$I_1$ of $\Gamma$ in $\Delta$ is $d_1$-dimensional, and that
the interpretation
$I_2$ of $\Delta$ in $\Gamma$ is $d_2$-dimensional.
Let $\phi(x,y_{1,1},\dots,y_{d_1,d_2})$ be the primitive positive formula
that defines $R_{I_2I_1}$ in $\Delta$.
Note that the primitive positive $\tau$-formula $\psi_{I_2}$ has $kd_2$ free variables; we can assume without loss of generality that $\psi_{I_2}$ only contains
unnested atomic formulas as conjuncts. Let
$(\psi_{I_2})_{I_1}$ be the primitive
positive $\sigma'$-formula obtained from $\psi_{I_2}$ by replacing each
conjunct $\psi'$ of $\psi_{I_2}$ by $(\psi')_{I_1}$,
and pushing existential quantifiers to the front. Then the formula
\begin{align*}
\exists y_{1,1}^1,\dots,y_{d_1,d_2}^k \; & \big ( \bigwedge_{i \leq k}
\phi(x_i,y^i_{1,1},\dots,y^i_{d_1,d_2})
\wedge \; (\psi_{I_2})_{I_1}
(y_{1,1}^1,\dots,y_{d_1,d_2}^1,\dots,y_{1,1}^k,\dots,y_{d_1,d_2}^k)
\big )
\end{align*}
is a primitive positive $\sigma'$-formula that defines
$\psi(x_1,\dots,x_k)$ over $\Delta'$.
\end{proof}

\section{Primitive Positive Interpretation Examples}
\label{sect:examples2}

\subsection*{Example 1}
Let $\Gamma$ be the structure with domain 
${\mathbb N}^2$ and a single binary relation $M := \{((u_1,u_2),(v_1,v_2)) \; | \; u_2=v_1 \text{ and } u_1,u_2,v_1,v_2 \in {\mathbb N} \}$.
Then $\Gamma$ and the structure $\Delta := ({\mathbb N};=)$ 
are primitive positive bi-interpretable. The interpretation $I$ of $\Gamma$ in $\Delta$ is 2-dimensional, the domain formula is \emph{true}, and the coordinate
map $h$ is the identity. The interpretation $J$ of $\Delta$ in $\Gamma$ is 1-dimensional, the domain formula is \emph{true}, and the coordinate map $g$
sends $(x,y)$ to $x$. Both interpretations are clearly primitive positive. 
Then $g(h(x,y))=z$ is definable by the formula $x=z$,
and $h(g(u),g(v))=w$ is primitive positive definable
by $$M(w,v) \wedge \exists p \; (M(u,p) \wedge M(w,p)) \; .$$  

\subsection*{Example 2}
An instructive example of two structures $\Gamma$ and $\Delta$ that are \emph{not} primitive
positive bi-interpretable, even though they are mutually primitive positive interpretable, is 
 $$\Gamma := \big({\mathbb N}^2; \{((u_1,u_2),(v_1,v_2)) \, | \, u_1 = v_1 \text{ and } u_1,u_2,v_1,v_2 \in {\mathbb N}\}\big)$$
 and $\Delta := ({\mathbb N};=)$. 
The two structures are not even \emph{first-order} 
bi-interpretable.
To see this, observe that the binary relation of $\Gamma$
is an equivalence relation, 
and that $\Aut(\Gamma)$ has 
a proper closed normal subgroup that is distinct from the one-element group, namely the set of all permutations that setwise fix 
the equivalence classes of this equivalence relation.
On the other hand,  $\Aut(\Delta)$ is the symmetric permutation group of a countably infinite set, which has no proper closed normal
subgroup that is distinct from the one-element group (it has exactly four proper normal subgroups~\cite{SchreierUlam}, of which only the one-element subgroup is closed).

\subsection*{Example 3}
The image of a continuous homomorphism $\xi$ 
from $\Pol(\Gamma)$ to $\Pol(\Delta)$ might be dense in 
$\Pol(\Delta)$ without being surjective, for $\omega$-categorical structures $\Gamma$ and $\Delta$. 
The basic idea of this example is due to Dugald Macpherson, 
and can be found in~\cite{HodgesLong} (on page~354).
Let $\Gamma$ be the structure $({\mathbb Q};<,P,P_4)$ where 
\begin{itemize}
\item $<$ is the usual strict order of the rational numbers,
\item $P \subseteq {\mathbb Q}$ is such that both
$P$ and $Q :={\mathbb Q} \setminus P$ are dense in $({\mathbb Q};<)$, and 
\item$P_4$ is the relation $\{(x_1,x_2,x_3,x_4) \in {\mathbb Q}^4 \; | \; x_1 = x_2 \text{ or } x_3 = x_4\}$.  
\end{itemize}
It is a well-known fact that all polymorphisms of $\Gamma$ are
essentially unary\footnote{A function $f \colon D^l \rightarrow D$ is called \emph{essentially unary} iff there exists an $i \leq l$ and a function $g \colon D \rightarrow D$ such that $f(x_1,\dots,x_l) = g(x_i)$ for all $x_1,\dots,x_l \in D$.} since they have to preserve $P_4$ (see 
e.g.~Lemma 5.3.2 in~\cite{Bodirsky-HDR}). 
The substructure $\Delta$ induced
by $P$ in $\Gamma$ has the primitive positive interpretation $(1,P,\id)$ in $\Gamma$. 
And indeed, since all functions of $\Pol(\Gamma)$ are essentially unary, the mapping which sends every unary function $f$ of $\Pol(\Gamma)$ to $f\rest_P$ induces a function $\xi$ from $\Pol(\Gamma)$
to $\Pol(\Delta)$ which is a continuous homomorphism and 
whose image is dense in $\Pol(\Delta)$. We claim that
$\xi$ is not surjective. 

A \emph{Dedekind cut $(S,T)$ of $P$} is a partition of $P$ into subsets $S,T$ with the property
that for all $s \in S$, $t \in T$ we have $s < t$. 
Those cuts are obtained by choosing either an irrational number
$r \in \mathbb R \setminus \mathbb Q$ or an element $r \in Q$, and setting $S := \{a \in P \; | \; a < r\}$ and $T := \{a \in P \; | \; a > r\}$.
Let $(S_1,T_1)$ be a Dedekind cut obtained from an element $q$ 
in $Q$, and let $(S_2,T_2)$ be a Dedekind cut obtained from an irrational number $i$. By a standard back-and-forth argument, there exists 
an $\alpha \in \Aut((P,<))$ that maps $S_1$ to $S_2$ and $T_1$ to $T_2$. Suppose for contradiction
that there is $\beta \in \Aut(\Gamma)$ with $\beta\rest_P= \alpha$. Then $s < \beta(q) < t$ for all $s \in S_2$, $t \in T_2$,
contradicting the irrationality of $i$.

\ignore{
\vspace{.4cm}
{\bf Example 4.} We adapt the previous example to demonstrate that 
$\Delta$ might have a primitive positive interpretation in an $\omega$-categorical structure $\Gamma$, but $\Delta$ is \emph{not}
the reduct of an $\omega$-categorical structure $\Delta'$ such
that there is a \emph{surjective} homomorphism from 
$\Pol(\Gamma)$ to $\Pol(\Delta')$.  

Let $\Gamma$ be the structure $({\mathbb Q};<,P,E)$ where
$<$ is the usual ordering of ${\mathbb Q}$, 
$P$ is a dense subset of ${\mathbb Q}$ such that ${\mathbb Q} \setminus P$ is also dense in ${\mathbb Q}$, and
$E$ is a subset of $P^2$ such that $(P,E)$ induces the random graph. Then
the structure $\Delta := (P;<,E)$ has the primitive positive interpretation $(1,P,\id)$ in $\Gamma$. Let
$\Delta'$ be any $\omega$-categorical structure such that
$\Delta$ is the reduct of $\Delta'$ and such that 
there is a continuous clone homomorphism $\xi$ from
$\Pol(\Gamma)$ to $\Pol(\Delta)$. We claim that 
$\xi$ cannot be surjective.  
}

\section{Constraint Satisfaction Problems}
\label{sect:csp}
Primitive positive interpretations play an important role in the study of the computational complexity of constraint satisfaction problems. 
For a structure 
$\Gamma$ with finite relational signature $\tau$, the \emph{constraint 
satisfaction problem for $\Gamma$} (denoted by $\Csp(\Gamma)$) is the 
computational problem to decide whether a given primitive positive 
$\tau$-sentence (that is, a primitive positive formula without free variables) is true in $\Gamma$. 
For example, when $\Gamma = (\{0,1,2\}; \neq)$, then
$\Csp(\Gamma)$ is the 3-coloring problem. 
When $\Gamma = ({\mathbb Q}; <)$, then $\Csp(\Gamma)$
is the acyclicity problem for finite directed graphs. 
Many computational problems studied in qualitative reasoning in artificial intelligence, but also in many other areas of theoretical computer science,  can be formulated as constraint satisfaction problems for $\omega$-categorical structures.

The subclass of problems of the form $\Csp(\Gamma)$ for finite 
$\Gamma$ attracted considerable interest in recent years. Feder and Vardi~\cite{FederVardi} conjectured that
such CSPs are either in P, or NP-complete. 
A very fruitful approach to this conjecture is the so-called
\emph{universal-algebraic approach}. 
One of the basic insights of this approach is that for finite $\Gamma$, the
complexity of $\Csp(\Gamma)$ only depends on the
pseudovariety generated by any of the polymorphism algebras of $\Gamma$. For $\omega$-categorical $\Gamma$, the same statement follows from Theorem~\ref{thm:pp-pseudovar} and
the following, which can be seen as a different formulation of results obtained in~\cite{JBK}.

\begin{theorem}[from~\cite{BodirskySurvey}]\label{thm:pp-interpret-reduce}
Let $\Gamma$ and $\Delta$ be structures with finite relational signatures.
If there is a primitive positive interpretation of $\Gamma$
in $\Delta$, then there is a polynomial-time reduction from
$\Csp(\Gamma)$ to $\Csp(\Delta)$. 
\end{theorem}

For finite structures $\Gamma$, 
this also shows that the complexity of $\Csp(\Gamma)$
is captured by the abstract polymorphism clone of $\Gamma$; see
Theorem~\ref{thm:birkhoff}. 
In other words, if $\Gamma$ and $\Delta$ are such that
their abstract polymorphism clones are isomorphic, then
$\Csp(\Gamma)$ and $\Csp(\Delta)$ are polynomial-time equivalent.  


Corollary~\ref{cor:csp-topo-clone}
gives a generalization of this fact for $\omega$-categorical structures: the complexity of $\Csp(\Gamma)$ only depends on the topological polymorphism clone of $\Gamma$.
In the following we explain that this is not only a fact of theoretical interest,
but that Theorem~\ref{thm:main} also provides a practical tool to prove hardness 
of $\Csp(\Gamma)$. An example will be given in Section~\ref{sect:examples3}. 
 
 Note that all algebras with domain of size at least two and with the property that all their operations are projections have, up to isomorphism,  the same abstract clone,
 which we denote by $\bf 1$.  For $1\leq i\leq k$, we denote the element of $\bf 1$ which correponds to the $k$-ary projection onto the $i$-th coordinate by
 $\pi_i^k$. So $\{\pi_i^k \; | \; i,k \in {\mathbb N}, i \leq k\}$ is the set of elements of the abstract clone $\bf 1$. Note that the topology on $\bf 1$ is the discrete topology since $\bf 1$ has only finitely many elements for each arity.

An example of a structure whose polymorphism clone is isomorphic to $\bf 1$ is the structure $(\{0,1\}; \OIT)$, where $\OIT := \{(0,0,1),(0,1,0),(1,0,0)\}$. The CSP for this structure is the well-known 
positive 1-IN-3-3SAT problem, which can be found in~\cite{GareyJohnson} and which is NP-complete.

\begin{theorem}\label{thm:find-1}
Let $\Gamma$ be an $\omega$-categorical structure. Then the following are equivalent.
\begin{enumerate}
\item All finite structures have a primitive positive interpretation
in $\Gamma$.
\item The structure $(\{0,1\}; \OIT)$ has a primitive positive interpretation in $\Gamma$.
\item $\Gamma$ has a polymorphism algebra $\C$ such that
the pseudovariety generated by $\C$ contains a two-element algebra $\A$ all of whose operations are projections.
\item There exists a continuous homomorphism from 
$\Pol(\Gamma)$ to ${\bf 1}$. 
\end{enumerate}
If one of those conditions applies, and $\Gamma$ has a relational signature, then $\Gamma$ has a finite signature
reduct $\Gamma'$ such that $\Csp(\Gamma')$ is NP-hard. 
\end{theorem}
\begin{proof}
The equivalence of $(1)$ and $(2)$ with (4) follows from Theorem~\ref{thm:main}, and the equivalence of $(3)$ with $(4)$ from Theorem~\ref{thm:topo-birk}. We remark that the equivalence between (1), (2) and $(3)$ can also be found in~\cite{Bodirsky-HDR}. 

To prove the statement about NP-hardness, let $\Gamma'$ be the reduct of
$\Gamma$ that contains exactly those relations that appear in the formulas
of the primitive positive interpretation of 
$(\{0,1\}; \OIT)$ in $\Gamma$. Note that $\Gamma'$ has finite signature, and still interprets 
$(\{0,1\}; \OIT)$ primitively positively. NP-hardness of $\Csp(\Gamma')$ now follows from the
mentioned fact that $\Csp((\{0,1\}; \OIT))$ is NP-hard, and Theorem~\ref{thm:pp-interpret-reduce}.
\end{proof}

\section{Constraint Satisfaction Example}
\label{sect:examples3}

Consider the structure $\Gamma = ({\mathbb Q}; \Betw)$ where $\Betw$ is the ternary relation $\{(x,y,z) \in {\mathbb Q}^3 \; | \; x<y<z \vee z<y<x\}$. Then $\Csp(\Gamma)$ is a well-known NP-complete problem known as the \emph{Betweenness problem}~\cite{Opatrny,GareyJohnson}. 
Applying the method presented in Section~\ref{sect:csp}, we will show NP-hardness of this problem by exhibiting a continuous clone homomorphism $\xi$ from $\Pol(\Gamma)$ to $\bf 1$.

In the following, for $k\geq 1$ and $x,y\in \Gamma^k$, we write ${{\neq}}(x,y)$ iff $x_j\neq y_j$ for all $1\leq j\leq k$.
\noindent {\bf Claim.} 
Let $k\geq 1$ and let $f\in\Pol(\Gamma)$ be $k$-ary. Then one of the following holds:

\begin{itemize}
\item[(1)] there is $1\leq d \leq k$  such that $f(x) < f(y)$ for all $x,y\in \Gamma^k$ with ${\neq}(x,y)$ and $x_d<y_d$;
\item[(2)] there is $1\leq d \leq k$  such that  $f(x) > f(y)$ for all $x,y\in \Gamma^k$ with ${\neq}(x,y)$ and $x_d<y_d$.
\end{itemize}
Since $d$ is clearly unique for each $f$,  
setting $\xi(f):=\pi^{k}_{d}$ defines a function $\xi$ from $\Pol(\Gamma)$ onto $\bf 1$. It is straightforward to check that $\xi$ is a homomorphism. To see that $\xi$ is continuous, observe that for $1\leq d\leq k$ the preimage of any $\pi^k_d$ under $\xi$ equals the intersection of $\Pol(\Gamma)$ with the set of all $k$-ary functions on $\Gamma$ which satisfy either (1) or (2); since the set of functions satisfying (1) or (2) is closed, so is $\xi\inv[\{\pi^k_i\}]$.

So we are left with the proof of the claim above. Observe first that either $f(0,\ldots,0)<f(1,\ldots,1)$ or $f(0,\ldots,0)>f(1,\ldots,1)$ holds: for if the two values were equal, this would contradict $\Betw(f(0,\ldots,0),f(1,\ldots,1),f(2,\ldots,2))$. We will now show that $f(0,\ldots,0)<f(1,\ldots,1)$ implies (1); then by symmetry of the statements, $f(0,\ldots,0)>f(1,\ldots,1)$ implies (2).
 
Observe the following: whenever $a,a',b,b'\in \Gamma^k$ are so that ${\neq}(a,a')$, ${\neq}(b,b')$, and $a_i<a_i'$ iff $b_i<b_i'$ for all $1\leq i\leq k$, then $f(a)<f(a')$ iff $f(b)<f(b')$. To see this, suppose without loss of generality that $f(a)<f(a')$ and $f(b)\geq f(b')$. Pick for all $1\leq i\leq k$ any $c_i<\min(a_i,a_i',b_i,b_i')$ if $a_i<a_i'$, and  $c_i>\max(a_i,a_i',b_i,b_i')$ otherwise. Pick moreover $d_i<\min(a_i,a_i',b_i,b_i')$ if $a_i>a_i'$, and  $d_i>\max(a_i,a_i',b_i,b_i')$ otherwise. Now $\Betw(c_i,a_i,a'_i)$ for all $1\leq i\leq k$ and $f(a)<f(a')$ imply $f(c)<f(a)$; likewise, $\Betw(a_i,a'_i,d_i)$  for all $1\leq i\leq k$ and $f(a)<f(a')$ imply $f(a')<f(d)$, and so $f(c)<f(d)$. However, the same argument with $b$ and $b'$ yields $f(d)<f(c)$, a contradiction. 

Now suppose that (1) does not hold, and let $c^0\in\Gamma^k$ be arbitrary. We will inductively define tuples $c^1,\ldots,c^k\in\Gamma^k$ such that $f(c^0)\geq f(c^1)\geq\cdots\geq f(c^k)$ and such that $c^0_i<c^k_i$ for all $1\leq i\leq k$, which contradicts our observation since $f(0,\ldots,0)<f(1,\ldots,1)$. For $0\leq j< k$, we define $c^{j+1}$ from $c^j$ as follows. Consider $x, y\in \Gamma^k$ witnessing the failure of (1) for $d=j$; that is, ${\neq}(x,y)$, $x_j<y_j$, and $f(x) \geq f(y)$ hold. Select $t\in \Gamma^k$ such that ${\neq}(c^j,t)$ and such that $c_i^j<t_i$ iff $x_i<y_i$ for all $1\leq i\leq k$. Then $c_j^j<t_j$, and the observation shows that $f(c^j)\geq f(t)$. For $1\leq i\leq k$, set $c^{j+1}_i:=c^j_i+k$ if $t_i> c^j_i$, and $c^{j+1}_i:=c^j_i-1$ otherwise. By our observation, $f(c^j)\geq f(c^{j+1})$; and since in the process every coordinate is increased by $k$ at least once, and decreased by $1$ at most $k-1$ times, we have $c^0_i<c^k_i$ for all $1\leq i\leq k$.

\section{Discussion}
\label{sect:conclusion}
Our results demonstrate that many properties of an $\omega$-categorical structure $\Gamma$ are already determined 
by the polymorphism clone of $\Gamma$ viewed as
a \emph{topological clone}, i.e., viewed as an abstract algebraic structure additionally equipped with the topology of pointwise convergence. 
One might ask which properties of
$\Gamma$ are captured by the abstract algebraic structure of the polymorphism clone
of $\Gamma$ \emph{without the topology}.  Observe that for finite $\Gamma$, the two concepts coincide.

We would like to point out that there is a considerable literature 
about $\omega$-categorical structures where the topology
on the automorphisms
is uniquely determined by the
abstract automorphism group; this is for instance the case if 
$\Aut(\Gamma)$ has the so-called \emph{small index property},
that is, all subgroups of countable index 
are open. (This is equivalent to saying that all homomorphisms
from $\Aut(\Gamma)$ to $S_\infty$, the symmetric group on a countably infinite set, are continuous.) 
The small index property has for instance been shown
\begin{itemize}
\item for $\Aut({\mathbb N};=)$
by Dixon, Neumann, and Thomas~\cite{DixonNeumannThomas};
\item for $\Aut({\mathbb Q};<)$ and for the automorphism group of the atomless Boolean algebra by Truss~\cite{Truss}; 
\item for the automorphism groups of the random graph~\cite{HodgesHodkinsonLascarShelah};
\item for all $\omega$-categorical $\omega$-stable structures~\cite{HodgesHodkinsonLascarShelah};
\item for the automorphism groups of the Henson graphs by Herwig~\cite{Herwig98}.
\end{itemize}
An example of two $\omega$-categorical structures (with infinite relational signature) whose automorphism groups are isomorphic as abstract groups but not as topological groups can be found in~\cite{EvansHewitt}. 

It is well-known that every Baire measurable  homomorphism between Polish groups is continuous (see e.g.~\cite{Kechris}).
So let us remark that there exists a model of ZF+DC where every set is Baire measurable~\cite{Shelah84}. 
For the 
structures $\Gamma$ that we need to model computational problems as $\Csp(\Gamma)$
it therefore seems fair to assume that the abstract automorphism
group of $\Gamma$ always determines the topological automorphism group (thanks to Todor Tsankov for pointing this out to us; consistency of this statement with ZF has already been observed in~\cite{Lascar}). 
However, this does not answer the question in which situations the abstract polymorphism clone determines the topological polymorphism clone.

In Theorem~\ref{thm:birkhoff}, if $\B\in  \HSPfin(\A)$, then $\B$ is in fact a homomorphic image of a subalgebra of $\A^n$, where $n=|\A|^{|\B|}$; that is, we have an explicit bound for the size of the power of $\A$ we have to take in order to represent $\B$. Peter Cameron has asked us whether a bound was known also in the locally oligomorphic case, i.e., in Theorem~\ref{thm:topo-birk}. By its nature of a compactness argument, our proof does not provide such a bound, and it would be interesting to find out whether a bound could be given also in our case.

\bibliographystyle{alpha}
\bibliography{topobirkhoff.bib}

\end{document}